\documentclass[12pt]{article}
\usepackage[margin=1in]{geometry}
\usepackage{fancyhdr}
\usepackage[skins, breakable]{tcolorbox}
\usepackage{amsmath, amsthm, amssymb}
\usepackage{amsfonts}
\usepackage{pgfplots}
\usepackage{tikz}
\usepackage[table,xcdraw]{xcolor}
\usepackage[parfill]{parskip}
\tcbuselibrary{breakable}
\usepackage{caption}
\usepackage{hyperref}
\hypersetup{
    colorlinks=true,
    linkcolor=blue,
    filecolor=magenta,      
    urlcolor=cyan,
    pdftitle={Overleaf Example},
    pdfpagemode=FullScreen,
    }

\newtheorem{theorem}{Theorem}[section]
\newtheorem{conjecture}[theorem]{Conjecture}
\newtheorem{corollary}[theorem]{Corollary}

\newtcolorbox{colorcorollary}[1]{colback=blue!5!white,
colframe=blue!75!black,fonttitle=\bfseries,
title={#1}}
\newtcolorbox{colortheorem}[1]{colback=black!5!white,
colframe=black!75!black,fonttitle=\bfseries,
title={#1}}
\newtcolorbox{colorconjecture}[1]{colback=orange!5!white,
colframe=orange!75!black,fonttitle=\bfseries,
title={#1}}
\newtcolorbox{colormain}[1]{colback=blue!5!white,
colframe=blue!75!black,fonttitle=\bfseries,
title={#1}}

\newcommand{\p}{\partial}

\usepackage[backend=biber,style=numeric]{biblatex}
\DeclareFieldFormat{labelnumber}{\thefield{entrykey}}
\begin{document}
\pagestyle{fancy}
\fancyhead{} 
\fancyhead[L]{26 September 2025}
\fancyhead[R]{McConkey, Seaton}
\fancyhead[C]{\textbf{Satellite Operations and $\theta$}}

\title{\textbf{Satellite Operations and $\theta$}}

\author{Rob McConkey \footnote{Department of Mathematics and Physics, Colorado State University Pueblo \\
2200 Bonforte Blvd Pueblo, CO 81001}
\and 
Luke J Seaton \footnote{Department of Mathematics, Michigan State University\\ 
619 Red Cedar Road East Lansing, MI 48824}}
\date{26 September 2025}

\maketitle

\begin{abstract}
We study the behavior of the knot invariant $\theta$ under satellite operations. First, we prove that $\theta$ is additive under connected sum. We then introduce a computational tool to generate $t$-twisted Whitehead doubles and apply it to explore the case of untwisted Whitehead doubles. We propose a conjecture describing the behavior of $\theta$ on untwisted Whitehead doubles and verify the conjecture for the first 2977 prime knots. The pair of invariants $\Theta = (\Delta,\theta)$ was introduced by Bar-Natan and van der Veen, where $\Delta$ is the Alexander polynomial. The invariant $\theta$ is easily computable and effective at distinguishing knots. Further exploration of satellite operations and $\theta$ is proposed to reveal new patterns among cables and general satellites.
\end{abstract}

{\footnotesize \textit{Keywords: }knot, Whitehead double, Alexander polynomial, satellite knot, knot invariant}\\
{\footnotesize Mathematics Subject Classification 2000: 57K14}

\section{Introduction}
\begin{figure}[h]
    \centering
    \includegraphics[width=0.5\linewidth]{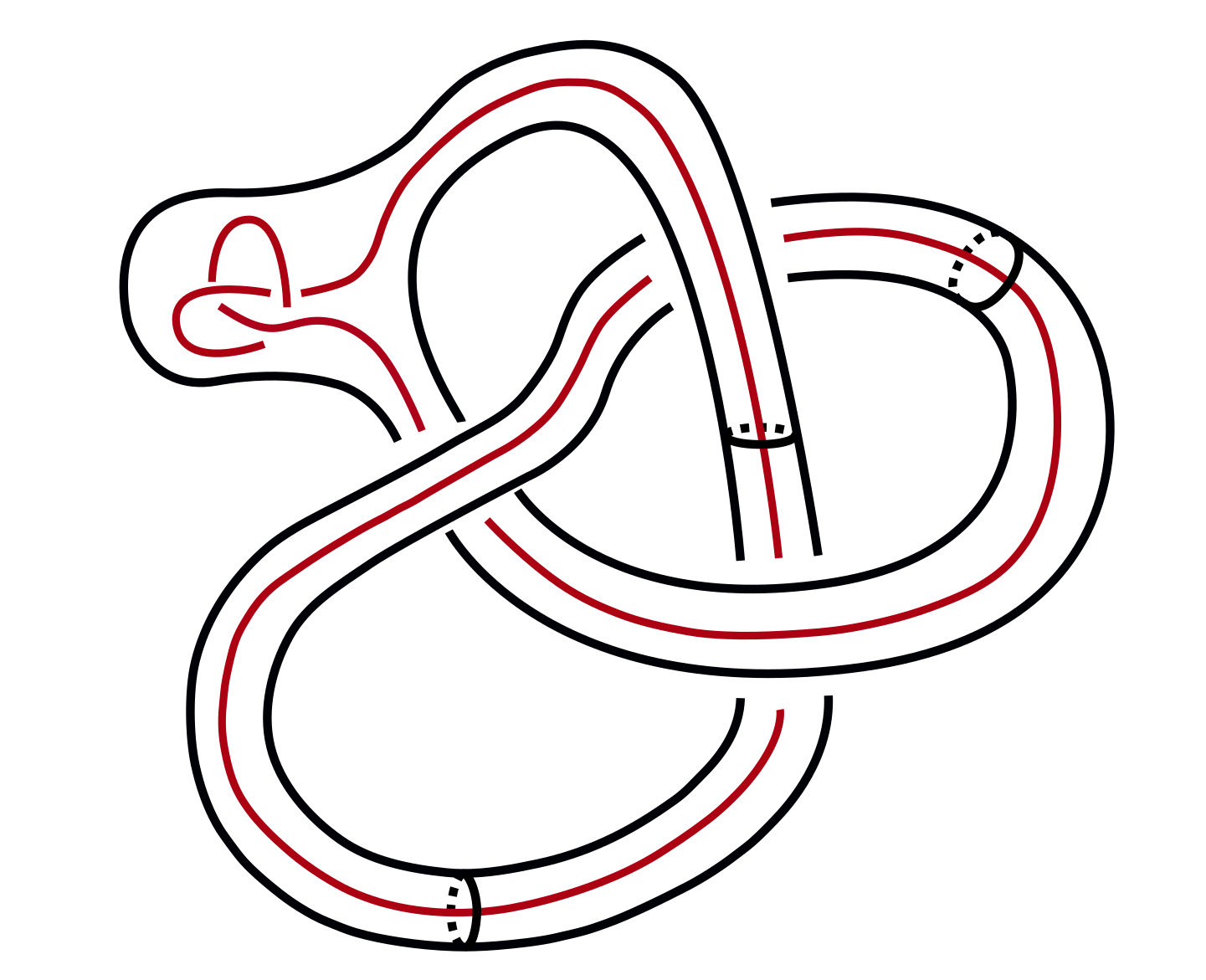}
    \caption{Connected sum of two trefoils as a satellite.}
    \label{fig:3_1and3_1}
\end{figure}

By work of Thurston \cite{Th82}, every knot in $S^3$ is either a torus knot, a hyperbolic knot, or a satellite knot. Satellite knots are frequently used as a testing ground for studying new knot invariants. Let the \textit{pattern knot} $P$ be a knot non-trivially wrapped inside a solid torus $V$, and $C$ be a classical knot called the \textit{companion knot}. The \textit{satellite knot} $C_P$ is the knot resulting from tying the solid torus $V$ into the shape of the companion $C$, as seen in Fig.\ref{fig:3_1and3_1}. Equivalently, $C_p$ is the result of gluing the boundary torus $\p V=T^2$ to the torus boundary of the exterior of the companion $\p (S^3\backslash\text{ nbd}(C))$. 


We define the \textit{wrapping number} to be the minimum number of times the pattern knot $P$ intersects a meridian disk of $V$. We say $P$ that is non-trivally wrapped inside the torus if it has wrapping number nonzero. We define the \textit{winding number} as the signed minimum intersection number where we take orientation into consideration. For example, the pattern knot in Fig. \ref{fig:3_1and3_1} has wrapping number one and winding number one. As a comparison, the pattern knot in Fig. \ref{fig:whiteheaddouble} has wrapping number two and winding number zero. 

Many invariants of satellite knots can be understood by examining their pattern and companion knots. For example, Theorem 6.15 in \cite{Li97} states a relationship for the Alexander polynomial shown in Eq. \ref{alexsatellite} where $r\geq 0$ is the linking number of $P$ with the meridian of the companion torus $\p V$. 

\begin{equation}\label{alexsatellite}
    \Delta_{C_P}(T) = \Delta_P(T)\Delta_C(T^r)
\end{equation}
The motivating question for this work asks: is there a satellite formula for $\theta$? 

The simplest satellite operation are those with wrapping number one. These operations are connected sum with a fixed knot $J$. In the connected sum, knots $K_1$ and $K_2$ are cut open, and the loose ends are joined so that the orientation is preserved in the sum, denoted $K_1\#K_2$. This connected sum with fixed $J=K_1$ (or $J=K_2$) can be viewed as a satellite operation. Either of $K_1$ or $K_2$ can be viewed as the pattern or companion where the pattern has a single strand wrap around the torus $V$. In Fig. \ref{fig:3_1and3_1}, the granny knot (a composite, satellite knot) is shown in red where both the pattern knot and the companion knot are right-handed trefoils. 

More interesting satellite operations arise when we increase the wrapping number of the pattern knot. When the pattern $P$ is a singly-clasped unknot with winding number 0 inside $V$, any companion $C$ yields a satellite knot called the Whitehead double of $C$, see Figure~\ref{fig:whiteheaddouble} for reference. Whitehead doubles are notable for having trivial Alexander polynomial. 

In the next sections, we study how a particular knot invariant $\theta$ behaves under satellite operations. First, we recall $\theta$ from work of Bar-Natan and van der Veen, and we show an example of how it can be computed in Section \ref{background}. In Section \ref{theta_additive}, we prove Theorem \ref{maintheorem} that states $\theta$ is additive under connected sum. We then extend our exploration of satellites to Whitehead doubles in Section \ref{WDs}. There we propose a conjecture describing the behavior of $\theta$ on untwisted Whitehead doubles and describe an algorithm to generate $t$-twisted Whitehead doubles. We verify the conjecture for the first 2977 prime knots (all knots with fewer than 13 crossings) using coded functions defined in Appendix \ref{appendix}. 

\section{Background on $\theta$}\label{background}
The pair of knot invariants $\Theta = (\Delta,\theta)$ was presented in The First International On-line Knot Theory Congress \cite{Ba25} by Bar-Natan, in joint work with van der Veen. In the process of writing this note, a preprint of their work was posted on the arXiv \cite{BV25}. With the familiar Alexander polynomial $\Delta$, the knot invariant $\theta$ is a strong candidate for future progress in knot theory. It is genuinely computable, as evident in our later computation of Whitehead doubles. It is conjectured to give an improved bound on genus, to be connected with other known invariants, and to potentially say something about ribbon knots. The pair $\Theta=(\Delta,\theta)$ can also distinguish more knots up to 15 crossings than the HOMFLY-PT polynomial and Khovanov homology combined (all while being more computable). To provide the background necessary for the arguments in Section \ref{theta_additive}, we recall the conventions for $\theta$ from \cite{Ba25} and provide an example computation. 

Consider an upright diagram $D$ of a knot $K$ where all of the crossings are oriented upwards and $D$ is a long knot diagram. We label each edge with a running index $k\in\{1, \dots , 2n+1\}$, where $n$ is the number of crossings in $D$, and a rotation number $\varphi_k$.\footnote{The signed number of full counterclockwise rotations of edge $k$ in the upright diagram. (For example, $\varphi_k=-1$ for the clockwise rotation of edge 2 in Fig. 3.)} Each crossing in $D$ is denoted by its sign and two lower edges:
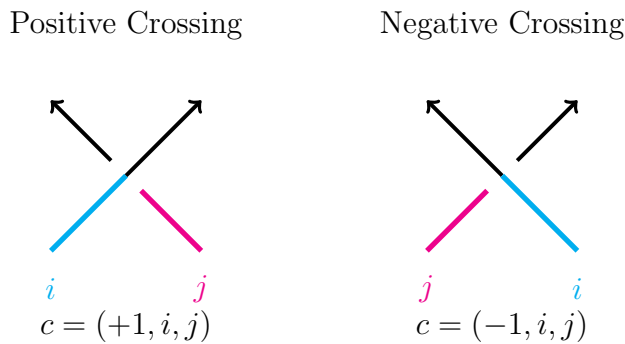
\begin{figure}[h]
    \begin{center}
    \begin{tikzpicture}
    \begin{scope} [xshift=5cm]
        \draw[line width=1.5pt, ->] (0,0) -- (-1,1); 
        \draw[line width=2pt,magenta] (-1,-1) -- (-0.2,-0.2); 
        \draw[line width=1.5pt, ->] (0.2,0.2) -- (1,1); 
        \draw[line width=2pt,cyan] (0,0) -- (1,-1); 
    \end{scope}
    
    \begin{scope}
        \draw[line width=1.5pt, ->] (-0.2,0.2) -- (-1,1); 
        \draw[line width=2pt,cyan] (-1,-1) -- (0,0); 
        \draw[line width=1.5pt, ->] (0,0) -- (1,1); 
        \draw[line width=2pt,magenta] (0.2,-0.2) -- (1,-1); 
    \end{scope}
    
    \node[cyan] at (-1,-1.5) {$i$};
    \node[magenta] at (1,-1.5) {$j$};
    \node[magenta] at (4,-1.5) {$j$};
    \node[cyan] at (6,-1.5) {$i$};
    \node at (0,2) {Positive Crossing};
    \node at (5,2) {Negative Crossing};
    \node at (0,-2) {$c=(+1,i,j)$};
    \node at (5,-2) {$c=(-1,i,j)$};    
\end{tikzpicture}
\end{center}
\caption{Crossings in an upright diagram.}
    \label{fig:uprightcrossings}
\end{figure}

We now introduce the traffic function $G = (g_{\alpha\beta})$. This $(2n+1)\times (2n+1)$ matrix records the traffic that enters at edge $\alpha$ and is measured at edge $\beta$. Cars travel through undercrossings with no problem, but at overcrossing bridges, they fall off with ``probability" $1-T^s$ (where $s$ is the sign of the crossing).
\newpage
\noindent \textbf{Example} We can compute the traffic function for a kink. There are three edges, so we have a $3\times 3$ matrix to compute.

\begin{center}
\begin{tikzpicture}
    \begin{scope} 
    \draw[line width=1.5pt] (0,-3) -- (0,-1); 
    \draw[line width=1.5pt] (0.2,-0.5) to[out=70,in=180] (1,0.5) to[out=0,in=90] (2,-0.5) 
        to[out=-90,in=0] (1,-1.5) to[out=180,in=-90] (0,-0.5); 
    \draw[line width=1.5pt, ->] (0,-0.5) -- (0,2); 
    \node at (1,-2) {$g_{11}=1$};
    \end{scope}
    \begin{scope}[scale=0.1,xshift=-2cm,yshift=-25cm,rotate=90]
    \fill[red] (-2,0) rectangle (2,1); 
    \fill[red] (-1,1) rectangle (1,1.8); 

    \fill[white] (-0.8,1) rectangle (-0.2,1.6); 
    \fill[white] (0.2,1) rectangle (0.8,1.6); 

    \fill[black] (-1.5,-0.3) circle (0.3); 
    \fill[black] (1.5,-0.3) circle (0.3); 

    \fill[gray] (-1.5,-0.3) circle (0.15);
    \fill[gray] (1.5,-0.3) circle (0.15);

    \fill[yellow] (2,0.5) circle (0.2); 
    \fill[gray] (-2,0.5) circle (0.2); 
    \end{scope}
    \begin{scope}[xshift=4cm]
    \draw[line width=1.5pt] (0,-3) -- (0,-1); 
    \draw[line width=1.5pt] (0.2,-0.5) to[out=70,in=180] (1,0.5) to[out=0,in=90] (2,-0.5) 
        to[out=-90,in=0] (1,-1.5) to[out=180,in=-90] (0,-0.5); 
    \draw[line width=1.5pt, ->] (0,-0.5) -- (0,2); 
    \node at (1,-2) {$g_{12}=T$};
    \end{scope}
    \begin{scope}[scale=0.1,xshift=38cm,yshift=-25cm,rotate=90]
    \fill[red] (-2,0) rectangle (2,1); 
    \fill[red] (-1,1) rectangle (1,1.8); 

    \fill[white] (-0.8,1) rectangle (-0.2,1.6); 
    \fill[white] (0.2,1) rectangle (0.8,1.6); 

    \fill[black] (-1.5,-0.3) circle (0.3); 
    \fill[black] (1.5,-0.3) circle (0.3); 

    \fill[gray] (-1.5,-0.3) circle (0.15);
    \fill[gray] (1.5,-0.3) circle (0.15);

    \fill[yellow] (2,0.5) circle (0.2); 
    \fill[gray] (-2,0.5) circle (0.2); 
    \end{scope}
    \begin{scope}[xshift=8cm]
    \draw[line width=1.5pt] (0,-3) -- (0,-1); 
    \draw[line width=1.5pt] (0.2,-0.5) to[out=70,in=180] (1,0.5) to[out=0,in=90] (2,-0.5) 
        to[out=-90,in=0] (1,-1.5) to[out=180,in=-90] (0,-0.5); 
    \draw[line width=1.5pt, ->] (0,-0.5) -- (0,2); 
    \node at (1,-2) {$g_{13}=1$};
    \end{scope}
    \begin{scope}[scale=0.1,xshift=78cm,yshift=-25cm,rotate=90]
    \fill[red] (-2,0) rectangle (2,1); 
    \fill[red] (-1,1) rectangle (1,1.8); 

    \fill[white] (-0.8,1) rectangle (-0.2,1.6); 
    \fill[white] (0.2,1) rectangle (0.8,1.6); 

    \fill[black] (-1.5,-0.3) circle (0.3); 
    \fill[black] (1.5,-0.3) circle (0.3); 

    \fill[gray] (-1.5,-0.3) circle (0.15);
    \fill[gray] (1.5,-0.3) circle (0.15);

    \fill[yellow] (2,0.5) circle (0.2); 
    \fill[gray] (-2,0.5) circle (0.2); 
    \end{scope}
    \begin{scope}[yshift=-6cm]
    \draw[line width=1.5pt] (0,-3) -- (0,-1); 
    \draw[line width=1.5pt] (0.2,-0.5) to[out=70,in=180] (1,0.5) to[out=0,in=90] (2,-0.5) 
        to[out=-90,in=0] (1,-1.5) to[out=180,in=-90] (0,-0.5); 
    \draw[line width=1.5pt, ->] (0,-0.5) -- (0,2); 
    \node at (1,-2) {$g_{21}=0$};
    \end{scope}
    \begin{scope}[scale=0.1,xshift=15cm,yshift=-55cm,rotate=-35]
    \fill[red] (-2,0) rectangle (2,1); 
    \fill[red] (-1,1) rectangle (1,1.8); 

    \fill[white] (-0.8,1) rectangle (-0.2,1.6); 
    \fill[white] (0.2,1) rectangle (0.8,1.6); 

    \fill[black] (-1.5,-0.3) circle (0.3); 
    \fill[black] (1.5,-0.3) circle (0.3); 

    \fill[gray] (-1.5,-0.3) circle (0.15);
    \fill[gray] (1.5,-0.3) circle (0.15);

    \fill[yellow] (2,0.5) circle (0.2); 
    \fill[gray] (-2,0.5) circle (0.2); 
    \end{scope}
    \begin{scope}[xshift=4cm,yshift=-6cm]
    \draw[line width=1.5pt] (0,-3) -- (0,-1); 
    \draw[line width=1.5pt] (0.2,-0.5) to[out=70,in=180] (1,0.5) to[out=0,in=90] (2,-0.5) 
        to[out=-90,in=0] (1,-1.5) to[out=180,in=-90] (0,-0.5); 
    \draw[line width=1.5pt, ->] (0,-0.5) -- (0,2); 
    \node at (1,-2) {$g_{22}=T$};
    \end{scope}
    \begin{scope}[scale=0.1,xshift=55cm,yshift=-55cm,rotate=-35]
    \fill[red] (-2,0) rectangle (2,1); 
    \fill[red] (-1,1) rectangle (1,1.8); 

    \fill[white] (-0.8,1) rectangle (-0.2,1.6); 
    \fill[white] (0.2,1) rectangle (0.8,1.6); 

    \fill[black] (-1.5,-0.3) circle (0.3); 
    \fill[black] (1.5,-0.3) circle (0.3); 

    \fill[gray] (-1.5,-0.3) circle (0.15);
    \fill[gray] (1.5,-0.3) circle (0.15);

    \fill[yellow] (2,0.5) circle (0.2); 
    \fill[gray] (-2,0.5) circle (0.2); 
    \end{scope}
    \begin{scope}[xshift=8cm,yshift=-6cm]
    \draw[line width=1.5pt] (0,-3) -- (0,-1); 
    \draw[line width=1.5pt] (0.2,-0.5) to[out=70,in=180] (1,0.5) to[out=0,in=90] (2,-0.5) 
        to[out=-90,in=0] (1,-1.5) to[out=180,in=-90] (0,-0.5); 
    \draw[line width=1.5pt, ->] (0,-0.5) -- (0,2); 
    \node at (1,-2) {$g_{31}=1$};
    \end{scope}
    \begin{scope}[scale=0.1,xshift=95cm,yshift=-55cm,rotate=-35]
    \fill[red] (-2,0) rectangle (2,1); 
    \fill[red] (-1,1) rectangle (1,1.8); 

    \fill[white] (-0.8,1) rectangle (-0.2,1.6); 
    \fill[white] (0.2,1) rectangle (0.8,1.6); 

    \fill[black] (-1.5,-0.3) circle (0.3); 
    \fill[black] (1.5,-0.3) circle (0.3); 

    \fill[gray] (-1.5,-0.3) circle (0.15);
    \fill[gray] (1.5,-0.3) circle (0.15);

    \fill[yellow] (2,0.5) circle (0.2); 
    \fill[gray] (-2,0.5) circle (0.2); 
    \end{scope}
    \begin{scope}[yshift=-12cm]
    \draw[line width=1.5pt] (0,-3) -- (0,-1); 
    \draw[line width=1.5pt] (0.2,-0.5) to[out=70,in=180] (1,0.5) to[out=0,in=90] (2,-0.5) 
        to[out=-90,in=0] (1,-1.5) to[out=180,in=-90] (0,-0.5); 
    \draw[line width=1.5pt, ->] (0,-0.5) -- (0,2); 
    \node at (1,-2) {$g_{31}=0$};
    \end{scope}
    \begin{scope}[scale=0.1,xshift=-2cm,yshift=-115cm,rotate=90]
    \fill[red] (-2,0) rectangle (2,1); 
    \fill[red] (-1,1) rectangle (1,1.8); 

    \fill[white] (-0.8,1) rectangle (-0.2,1.6); 
    \fill[white] (0.2,1) rectangle (0.8,1.6); 

    \fill[black] (-1.5,-0.3) circle (0.3); 
    \fill[black] (1.5,-0.3) circle (0.3); 

    \fill[gray] (-1.5,-0.3) circle (0.15);
    \fill[gray] (1.5,-0.3) circle (0.15);

    \fill[yellow] (2,0.5) circle (0.2); 
    \fill[gray] (-2,0.5) circle (0.2); 
    \end{scope}
    \begin{scope}[xshift=4cm,yshift=-12cm]
    \draw[line width=1.5pt] (0,-3) -- (0,-1); 
    \draw[line width=1.5pt] (0.2,-0.5) to[out=70,in=180] (1,0.5) to[out=0,in=90] (2,-0.5) 
        to[out=-90,in=0] (1,-1.5) to[out=180,in=-90] (0,-0.5); 
    \draw[line width=1.5pt, ->] (0,-0.5) -- (0,2); 
    \node at (1,-2) {$g_{32}=0$};
    \end{scope}
    \begin{scope}[scale=0.1,xshift=38cm,yshift=-115cm,rotate=90]
    \fill[red] (-2,0) rectangle (2,1); 
    \fill[red] (-1,1) rectangle (1,1.8); 

    \fill[white] (-0.8,1) rectangle (-0.2,1.6); 
    \fill[white] (0.2,1) rectangle (0.8,1.6); 

    \fill[black] (-1.5,-0.3) circle (0.3); 
    \fill[black] (1.5,-0.3) circle (0.3); 

    \fill[gray] (-1.5,-0.3) circle (0.15);
    \fill[gray] (1.5,-0.3) circle (0.15);

    \fill[yellow] (2,0.5) circle (0.2); 
    \fill[gray] (-2,0.5) circle (0.2); 
    \end{scope}
    \begin{scope}[xshift=8cm,yshift=-12cm]
    \draw[line width=1.5pt] (0,-3) -- (0,-1); 
    \draw[line width=1.5pt] (0.2,-0.5) to[out=70,in=180] (1,0.5) to[out=0,in=90] (2,-0.5) 
        to[out=-90,in=0] (1,-1.5) to[out=180,in=-90] (0,-0.5); 
    \draw[line width=1.5pt, ->] (0,-0.5) -- (0,2); 
    \node at (1,-2) {$g_{33}=1$};
    \end{scope}
    \begin{scope}[scale=0.1,xshift=78cm,yshift=-115cm,rotate=90]
    \fill[red] (-2,0) rectangle (2,1); 
    \fill[red] (-1,1) rectangle (1,1.8); 

    \fill[white] (-0.8,1) rectangle (-0.2,1.6); 
    \fill[white] (0.2,1) rectangle (0.8,1.6); 

    \fill[black] (-1.5,-0.3) circle (0.3); 
    \fill[black] (1.5,-0.3) circle (0.3); 

    \fill[gray] (-1.5,-0.3) circle (0.15);
    \fill[gray] (1.5,-0.3) circle (0.15);

    \fill[yellow] (2,0.5) circle (0.2); 
    \fill[gray] (-2,0.5) circle (0.2); 
    \end{scope}

    \draw[line width=0.5mm,blue] (-0.2,-1.5) -- (0.2,-1.5);
    \draw[line width=0.5mm,blue] (-0.2,-7.5) -- (0.2,-7.5);
    \draw[line width=0.5mm,blue] (-0.2,-13.5) -- (0.2,-13.5);
    \draw[line width=0.5mm,blue] (5.8,-0.8) -- (6.2,-0.9);
    \draw[line width=0.5mm,blue] (5.8,-6.8) -- (6.2,-6.9);
    \draw[line width=0.5mm,blue] (5.8,-12.8) -- (6.2,-12.9);
    \draw[line width=0.5mm,blue] (7.8,1) -- (8.2,1);
    \draw[line width=0.5mm,blue] (7.8,-5) -- (8.2,-5);
    \draw[line width=0.5mm,blue] (7.8,-11) -- (8.2,-11);

    \node at (11,-10) {\Large{$G = \begin{pmatrix}
    1 & T & 1\\ 0 & T & 1\\ 0 & 0 & 1
\end{pmatrix}$}};
\end{tikzpicture}
\end{center}
\begin{center}{\footnotesize Fig. 3. $\quad$ Calculating the traffic function for a kink.}\end{center}
\stepcounter{figure}
Note there are some cases where the traffic counter never sees any traffic (where $g_{\alpha\beta}=0$). In some places, such as $g_{12}$ and $g_{22}$ in this example, cars can get counted multiple times. If traffic enters at edge $\alpha =1$ and is counted on edge $\beta=2$, then all of the cars go under the bridge, all cars hit the counter, but some cars fall off the bridge and hit the counter again (and some again...). The traffic counter sees the following traffic: $$g_{12} = 1 + (1-T^{-1}) + (1-T^{-1})^2 + \cdots = \sum_{p=0}^\infty (1-T^{-1})^p = \frac{1}{1-(1-T^{-1})} = T$$

We will use the entries of the matrix $G$ to compute functions $F_1,F_2,$ and $F_3$ necessary for calculating $\theta$. Let $T_1$ and $T_2$ be indeterminates and set $T_3 := T_1T_2$. For $\nu=1,2,3$, let $G_\nu$ be the traffic function where $T$ is replaced by $T_\nu$. That is, the notation $g_{\nu \alpha\beta}$ refers to the traffic that enters at edge $\alpha$ and measured at edge $\beta$ in terms of $T_\nu$. Given crossings $c=(s,i,j),c_0=(s_0,i_0,j_0)$, and $c_1=(s_1,i_1,j_1)$, we have the functions
\begin{align*}
    F_1(c) &= s[1/2 - g_{3ii} +T_2^s g_{1ii}g_{2ji} - T^s_2g_{3jj}g_{2ji} - (T^s_2 - 1)g_{3ii}g_{2ji}\\
    &\qquad + (T^s_3 -1) g_{2ji}g_{3ji} - g_{1ii}g_{2jj} + 2g_{3ii}g_{2jj} + g_{1ii}g_{3jj} - g_{2ii}g_{3jj}]\\ 
    &\qquad + \frac{s}{T^s_2 -1}[(T^s_1-1)T^s_2 (g_{3jj}g_{1ji} - g_{2jj}g_{1ji} + T^s_2 g_{1ji}g_{2ji})\\
    &\qquad +(T^s_3-1)(g_{3ji} - T^s_2
g_{1ii}g_{3ji} + g_{2ij}g_{3ji} + (T^s_2 - 2)g_{2jj}g_{3ji})\\
&\qquad -(T^s_1-1)(T^s_2+1)(T^s_3-1)g_{1ji}g_{3ji}]\\
F_2(c_0,c_1) &= \frac{s_1(T^{s_0}_1-1)(T^{s_1}_3-1)g_{1j_1i_0}g_{3j_0i_1}}{T^{s_1}_2-1}(T^{s_0}_2g_{2i_1i_0} + g_{2j_1j_0} - T^{s_0}_2g_{2j_1i_0} - g_{2i_1j_0})\\
F_3(\varphi_k,k) &= \varphi_k(g_{3kk} - 1/2)
\end{align*}

\begin{colortheorem}
    {\begin{theorem}{[Main Theorem in \cite{Ba25}]}\end{theorem}}
    The following is a knot invariant: $$\theta(D) = \Delta_1\Delta_2\Delta_3\left( \sum_c F_1(c) + \sum_{c_0,c_1} F_2(c_0,c_1) + \sum_k F_3(\varphi_k, k) \right).$$
\end{colortheorem}

Note that $\Delta_1$ is the Alexander polynomial $\Delta$ of the knot $K$, where $D$ is a diagram of $K$, with indeterminate $T_1$. Similarly, $\Delta_2$ is in terms of $T_2$, and $\Delta_3$ is in terms of $T_3 = T_1T_2$. Later we write $N(K)$ to be a shorthand for this product of Alexander polynomials.

\begin{figure}[h]
    \centering
    \includegraphics[scale=0.25]{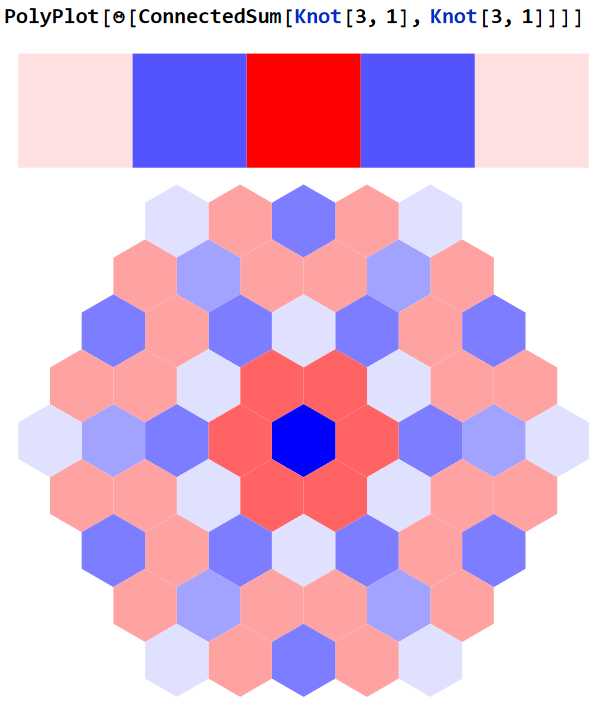}
    \includegraphics[scale=0.25]{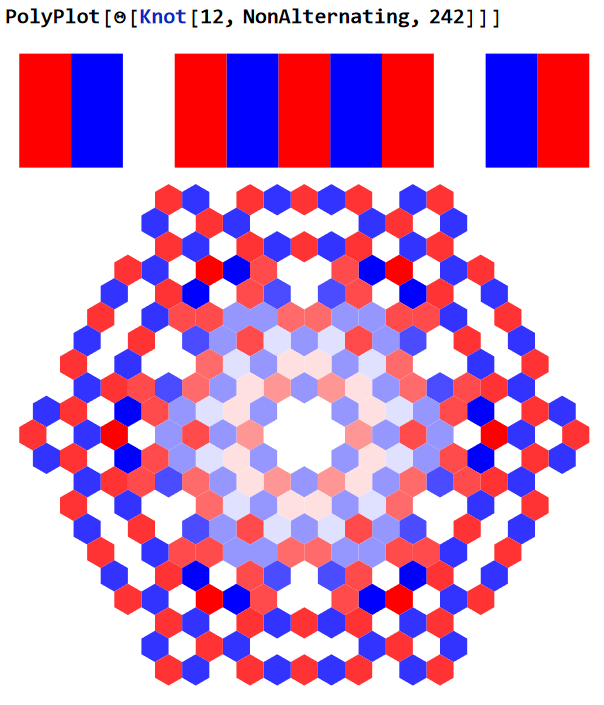}
    \includegraphics[scale=0.25]{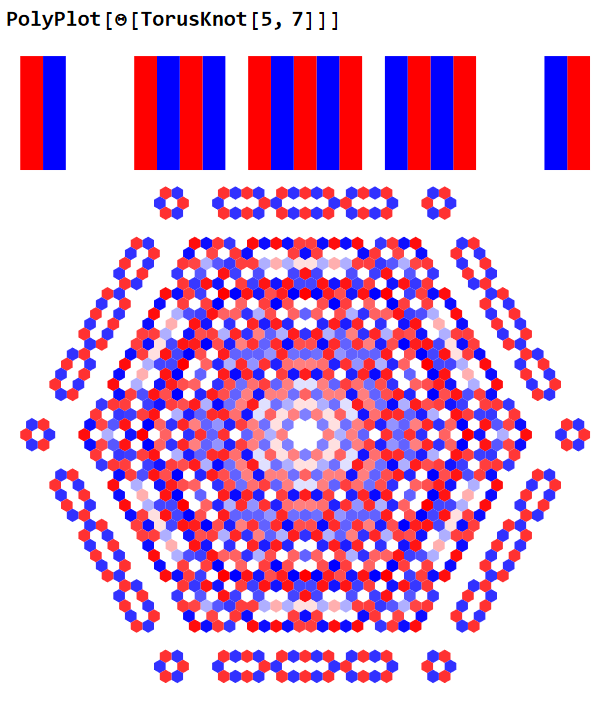}
    \caption{From left to right: $\Theta(3_1\#3_1)$, $\Theta(12n_242)$, and $\Theta(T_{5,7})$}
    \label{fig:thetaQR}
\end{figure}

\newpage
One of the major draws of this knot invariant is its beauty when plotted in powers of $T_1$ and $T_2$. The red or blue color is associated with positive or negative coefficients, respectively. The intensity of the color is proportional to the absolute value of the coefficient (i.e. brighter red means higher positive number, and white means zero). In Figure \ref{fig:thetaQR}, the bar code at the top represents the symmetric Alexander polynomial while the QR code represents $\theta$. 

\section{Connected Sum and $\theta$}\label{theta_additive}
The knot invariant $\rho_1$ (also introduced by Bar-Natan and van der Veen \cite{BV22}) uses the traffic function $G=(g_{\alpha\beta})$ and has been studied since 2022. Robert Quarles \cite{Qu22} proved a formula for $\rho_1$ of the connected sum of two knots in his PhD thesis. 

\begin{colortheorem}
    {\begin{theorem}{[Quarles, Theorem 3.2.2 in \cite{Qu22}]}\end{theorem}}
    For a knot $K_1\#K_2$, we have $$\rho_1 (K_1\#K_2) = \Delta^2_{K_2} \rho_1 (K_1) + \Delta^2_{K_1} \rho_1 (K_2).$$
\end{colortheorem}

Since $\rho_1$ is constructed similarly to $\theta$, the statement of the theorem inspired the idea of how to extend it to $\theta$. 

\begin{colormain}
{\begin{theorem}\label{maintheorem}\end{theorem}}
    Let $N(K) = \Delta_1(K)\Delta_2(K)\Delta_3(K)$, the normalization factor in the definition of $\theta$. For a knot $K_1\#K_2$, we have $$\theta (K_1\#K_2) = N(K_2)\theta(K_1) + N(K_1)\theta(K_2).$$
\end{colormain}

\begin{proof}
    We start with an upright diagram for a knot $K_1$ with $n_1$ crossings and stack an upright diagram for a knot $K_2$ with $n_2$ crossings on top of it. Figure \ref{fig:uprightconnectedsum} gives a diagram for $K_1\# K_2$. Let $\ell$ be the edge denoted $2n_1+1$ that connects $K_1$ to $K_2$. We now have an upright diagram with $2n_1+2n_2+1$ edges. 
    
    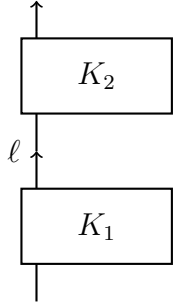
\begin{figure}[h]
    \begin{center}
    \begin{tikzpicture}
    \draw[thick] (0,0) rectangle (2,1) node at (1,0.5) {$K_1$};
    \draw[thick,->] (0.2,1) -- (0.2,1.5); 
    \draw[thick] (0.2,0) -- (0.2,-0.5); 
    
    \draw[thick] (0,2) rectangle (2,3) node at (1,2.5) {$K_2$};
    \draw[thick,->] (0.2,3) -- (0.2,3.5); 
    \draw[thick] (0.2,2) -- (0.2,1.5); 

    \node at (-0.1,1.5) {$\ell$};
    \end{tikzpicture}
    \end{center}
       \caption{An upright diagram for the connected sum of $K_1$ and $K_2$.}
        \label{fig:uprightconnectedsum}
    \end{figure}
    
    It suffices to show that the summations of $F_1, F_2, F_3$ break into two summations (one for $K_1$ and the other for $K_2$). We do this in three claims.

    \textbf{Claim 1:}
    \begin{align}\label{f1}
        \sum_{c\in K_1\#K_2} F_1(c) &= \sum_{c\in K_1} F_1(c) + \sum_{c\in K_2} F_1(c)
    \end{align}
    First, we consider the summation on the left side of Equation \ref{f1} over all crossings $c$ in $K_1\#K_2$. Each crossing $c$ exists in either $K_1$ or $K_2$, with no lower edges in common, so the $F_1$ summation splits as above.  
    
    \textbf{Claim 2:} 
    \begin{align}\label{f2}
        \sum_{c_0,c_1\in K_1\#K_2} F_2(c_0,c_1) = \sum_{c_0,c_1\in K_1} F_2(c_0,c_1) + \sum_{c_0,c_1\in K_2} F_2(c_0,c_1)
    \end{align}
    The function $F_2$ considers pairs of crossings. When both crossings are in $K_1$ (resp. both in $K_2$), the term appears in the first summation (resp. second summation) on the right side of Equation \ref{f2}. The case that remains is when one crossing is from $K_1$ and one is from $K_2$.
    
    We want to show that $F_2(c_1,c_2)=0$ when $c_1 = (s_1,i_1,j_1)$ is a crossing in $K_1$ and $c_2 = (s_2,i_2,j_2)$ is a crossing in $K_2$. Consider $F_2$:
    $$F_2(c_1,c_2) = \frac{s_2(T^{s_1}_1-1)(T^{s_2}_3-1)\textcolor{magenta}{g_{1j_2i_1}}g_{3j_1i_2}}{T^{s_2}_2-1}(T^{s_1}_2g_{2i_2i_1} + g_{2j_2j_1} - T^{s_1}_2g_{2j_2i_1} - g_{2i_2j_1})$$
    The factor $g_{1j_2i_1}$ (highlighted in \textcolor{magenta}{magenta}) above is the traffic measured at $i_1$ from inputting traffic at edge $j_2$. However, traffic input at an edge in $K_2$ will never travel down to any edge in $K_1$. Thus, $g_{1j_2i_1}=0$, so $F_2(c_1,c_2)=0$. Then we have equality above. 
    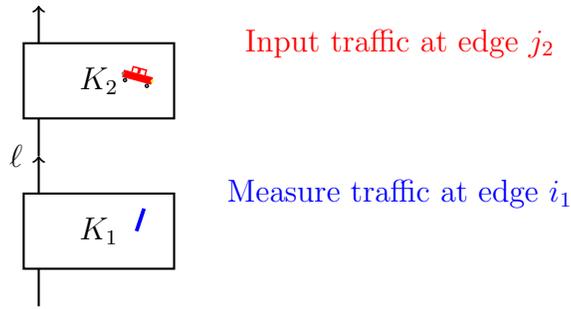
\begin{figure}[h]
    \begin{center}
    \begin{tikzpicture}
    \draw[thick] (0,0) rectangle (2,1) node at (1,0.5) {$K_1$};
    \draw[thick,->] (0.2,1) -- (0.2,1.5); 
    \draw[thick] (0.2,0) -- (0.2,-0.5); 
    
    \draw[thick] (0,2) rectangle (2,3) node at (1,2.5) {$K_2$};
    \draw[thick,->] (0.2,3) -- (0.2,3.5); 
    \draw[thick] (0.2,2) -- (0.2,1.5); 

    \node at (-0.1,1.5) {$\ell$};

    \draw[line width=0.5mm,blue] (1.5, 0.5) -- (1.6,0.8);

    \node[red] at (5,3) {Input traffic at edge $j_2$};
    \node[blue] at (5,1) {Measure traffic at edge $i_1$};

    \begin{scope}[scale=0.1,xshift=15cm,yshift=25cm,rotate=-15]
    \fill[red] (-2,0) rectangle (2,1); 
    \fill[red] (-1,1) rectangle (1,1.8); 

    \fill[white] (-0.8,1) rectangle (-0.2,1.6); 
    \fill[white] (0.2,1) rectangle (0.8,1.6); 

    \fill[black] (-1.5,-0.3) circle (0.3); 
    \fill[black] (1.5,-0.3) circle (0.3); 

    \fill[gray] (-1.5,-0.3) circle (0.15);
    \fill[gray] (1.5,-0.3) circle (0.15);

    \fill[yellow] (2,0.5) circle (0.2); 
    \fill[gray] (-2,0.5) circle (0.2); 
    \end{scope}
    \end{tikzpicture}
    \end{center}
    \caption{Traffic that enters at an edge in $K_2$ and is measured at an edge in $K_1$.}
        \label{fig:claim2}
    \end{figure}
    
\newpage
    \textbf{Claim 3: }
    \begin{align}
        \sum_{k\in K_1\#K_2} F_3(\varphi_k, k) &= \sum_{k \in K_1} F_3(\varphi_k, k) + \sum_{k \in K_2} F_3(\varphi_k, k)
    \end{align}
    Lastly, we consider $F_3$ which considers edges. Each edge is in either $K_1$ or $K_2$ except for edge $\ell$. This edge $\ell$ has rotation number $\varphi_\ell =0$, so $F_3(\ell) = 0$. So we have equality. 

    Now we can sum all three of the claims, and we get an equation:
    \begin{align}\label{sums}
        \sum_{K_1\#K_2} F_1 + \sum_{K_1\#K_2} F_2 + \sum_{K_1\#K_2} F_3 =& \left( \sum_{K_1} F_1 + \sum_{K_1} F_2 + \sum_{K_1} F_3 \right)\\
        & \qquad + \left( \sum_{K_2} F_1 + \sum_{K_2} F_2 + \sum_{K_2} F_3 \right) \notag
    \end{align}
    Next, we consider the definition of $\theta$. The Alexander polynomial is multiplicative with connected sum, so we have 
    \begin{align*}
        N(K_1\#K_2) &= \Delta_1(K_1\#K_2)\cdot\Delta_2(K_1\#K_2)\cdot\Delta_3(K_1\#K_2)\\
        &= [\Delta_1(K_1)\cdot\Delta_2(K_1)\cdot\Delta_3(K_1)]\cdot[\Delta_1(K_2)\cdot\Delta_2(K_2)\cdot\Delta_3(K_2)]\\
        &= N(K_1)N(K_2)
    \end{align*} 
    If we multiply the equation of sums (Eq. \ref{sums}) by $N(K_1\#K_2)$ on the left and $N(K_1)N(K_2)$ on the right, we arrive at our desired formula.
    $$\theta (K_1\#K_2) = N(K_2)\theta(K_1) + N(K_1)\theta(K_2)$$
\end{proof}

One may be tempted to apply the same argument to give an alternative proof for Quarles' theorem; however, $\rho_1$ uses upper edges from the crossings, such as in $g_{j,j^+}$, where $\theta$ uses only the bottom edges $i$ and $j$. The argument above does not immediately extend to show that $\rho_1$ is additive under connected sum.\\

\noindent \textbf{Note:} It is conjectured in \cite{Ba25} by Bar-Natan and van der Veen that $\theta (\overline{K}) = - \theta (K)$. While a direct argument for this depends on the conjectured hexagonal symmetry of $\theta$, it may be possible to prove using the follow result. 

\begin{colorcorollary}
    {\begin{corollary}\end{corollary}}
    Let $\overline{K}$ be the mirror of a knot $K$. $$\theta(K\#\overline{K}) =0 \iff \theta (\overline{K}) = - \theta (K)$$
\end{colorcorollary}

\begin{proof}
    ($\Rightarrow$) By the connected sum formula, $$0 = \theta (K\#\overline{K}) = N(\overline{K}) \theta(K) + N(K) \theta(\overline{K})$$
    Since $K$ and $\overline{K}$ have the same Alexander polynomial (Proposition 6.12 in \cite{Li97}), we have that $N(K) = N(\overline{K})$. Then $$0=N(K)(\theta(K) + \theta(\overline{K})).$$ Since $\Delta_K(1)=\pm 1$ for any knot $K$ (Theorem 6.10 in \cite{Li97}), we know that we always have $N(K)\neq 0$. Then we must have $\theta (\overline{K}) = - \theta (K)$.\\
    ($\Leftarrow$) Assume $\theta (\overline{K}) = - \theta (K)$ and the same arguments apply: 
    \begin{align*}
        \theta(K\#\overline{K}) &= N(\overline{K}) \theta(K) + N(K) \theta(\overline{K})\\
        &= N(\overline{K}) \theta(K) - N(K) \theta(K)\\
        &= \theta (K) (N(\overline{K}) - N(K))
    \end{align*}
    Since $N(\overline{K})=N(K)$ as before, we have $\theta(K\#\overline{K}) = 0$. 
\end{proof}

\section{Whitehead Doubles}\label{WDs}
\begin{figure}[h]
\centering
\begin{tikzpicture}
\begin{scope}
    \node at (0,0) {\includegraphics[width=0.6\linewidth]{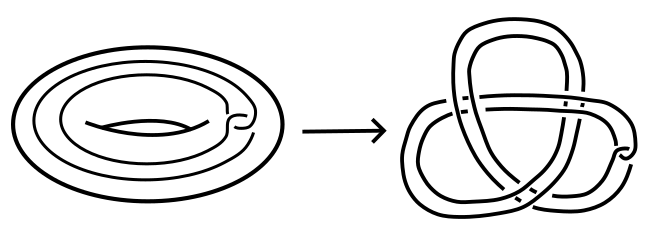}};
    \node at (-4.2,-0.6) {$P\subset V$};
    \node[draw,fill=white] at (-2.1,-0.6) {$t$};
    \node[draw,fill=white] at (2.2,1) {$t+3$};
\end{scope}
\end{tikzpicture}
\caption{The positively-clasped, $t$-twisted Whitehead double of the left-handed trefoil.}
\label{fig:whiteheaddouble}
\end{figure}

Another well-studied family of satellites is the Whitehead doubles. The pattern knot for a Whitehead double is the trivial knot wrapped around the torus with the ends clasped (as in the left of Fig.~\ref{fig:whiteheaddouble}. When tying (or gluing) the torus to create the satellite, twists are introduced to counteract the writhe $Wr$ of the projection. The $t$-twisted Whitehead double is created by adding $t$ full twists to the pattern that yields $t-Wr$ full twists in the Whitehead double. For example, the left-handed trefoil has writhe $-3$, so the $t$-twisted Whitehead double of the left-handed trefoil has $t+3$ twists in Fig.~\ref{fig:whiteheaddouble}. 

Note the relationship between positively- and negatively-clasped Whitehead doubles
\begin{equation}
    \overline{WD(+,\overline{K},-t)} = WD(-,K,t)
\end{equation}
where an overline represents the mirror. Using Eq. \ref{alexsatellite}, we obtain Alexander polynomial for the $t$-twisted Whitehead double \cite{Ro76}.
\begin{equation}
    \Delta_{WD(+,K,t)}(T) = -t\cdot T + (2t+1) -t\cdot T^{-1}
\end{equation}
In particular, untwisted Whitehead doubles have trivial Alexander polynomial. 

\subsection{Whitehead Doubles and $\theta$}
Inspired by Conjecture 4 from \cite{Ba25} that $\theta$ could be equal to Ohtsuki's ``two-loop polynomial'', we compared Ohtsuki's cabling formula (Theorem 4.1 in \cite{Oh04}) with examples of cabling for $\Theta$ in Mathematica. After not finding direct equality, we turned to a family of simpler satellites: Whitehead doubles. Based on the relationship between the cabling formula and the general satellite formula for the Alexander polynomial, we suspect that there is a similar relationship for $\Theta$. 

\begin{figure}
\begin{center}
    \includegraphics[scale=0.25]{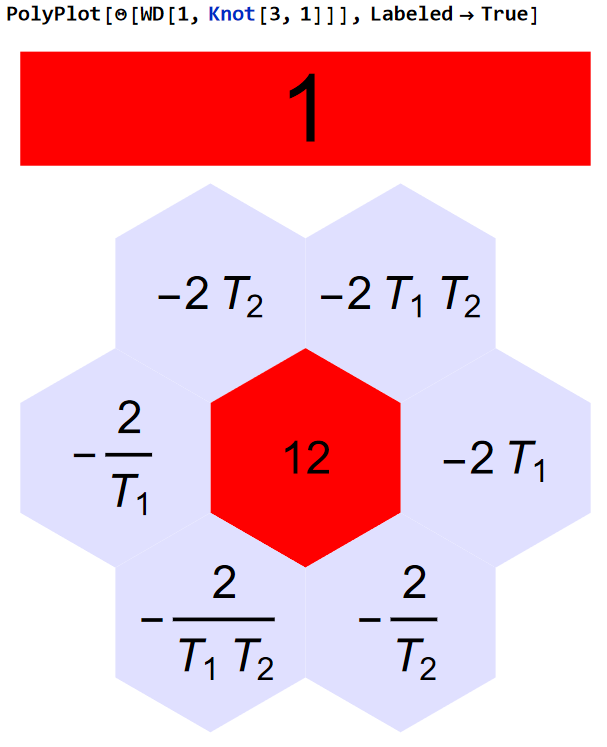}
    \includegraphics[scale=0.25]{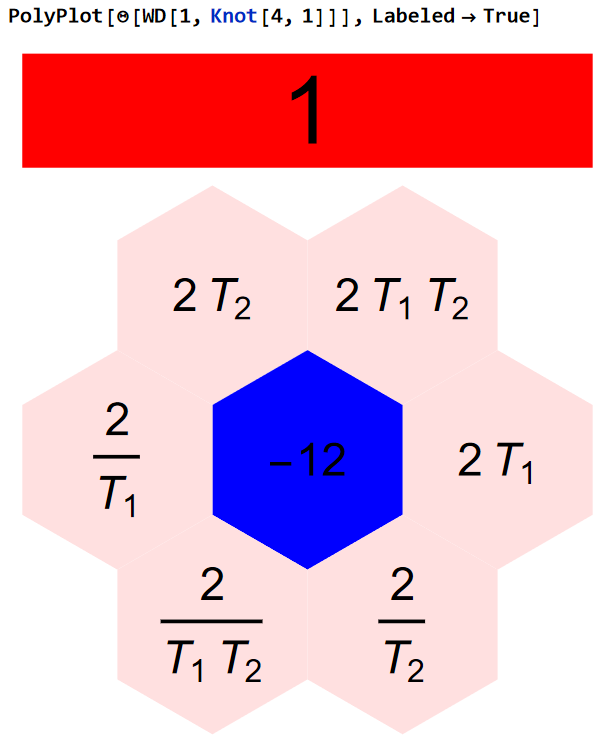}
    \includegraphics[scale=0.25]{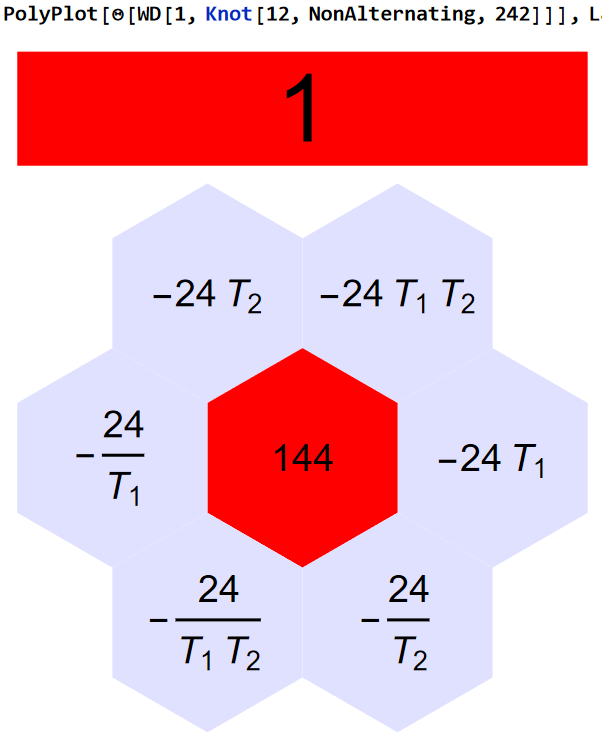}
\end{center}
\caption{$\Theta$ for Whitehead doubles of the trefoil, figure-8, and 12n\_242 knots.}
    \label{fig:flowers}
\end{figure}

We considered Whitehead doubles in an effort to simplify potential terms in a general satellite formula for $\Theta$. These satellites have the unknot as the pattern knot, so using Equation \ref{alexsatellite} we have \begin{equation}
    \Delta_{WD(K)}(T) = \Delta_U(T)\Delta_K(T^r) = \Delta_K(1)
\end{equation} 
because the linking number $r=0$. This fact further implies $N(WD(K))=1$ and $\Delta'(U) = \Delta'(WD(K)) = 0$. Ohtsuki's cabling formula includes several terms, some of which depend on $\Delta'_K$. However, $\Delta'_K(1)=0$ for all knots by symmetry of the Alexander polynomial and the Chain Rule: 
\begin{equation}
    \frac{d}{dt}\Delta_K(t) = \frac{d}{dt} \Delta(t^{-1}) = \Delta'_K (t^{-1})\cdot \frac{d}{dt}(t^{-1}) = -t^{-2} \Delta'_K(t^{-1})
\end{equation}
At $t=1$, we have 
\begin{equation}
    \Delta'_K(1) = -1^{-2}\Delta'_K(1^{-1}) = -\Delta'_K \quad\Longrightarrow\quad \Delta'_K(1)=0
\end{equation}
This led us to consider the second derivative (where $\Delta''_K(1)$ is often nonzero), and we found a pattern. 
\begin{colorconjecture}
{\begin{conjecture}\end{conjecture}}
    Let $K$ be any knot. Then $\theta$ of the Whitehead double of $K$ is 
    $$\theta (WD(+,K)) = q \left( 6-T_1-T_2-T_1T_2-\frac{1}{T_1} - \frac{1}{T_2} -\frac{1}{T_1T_2}\right)$$
    where the constant $q=\Delta_K''(1)$. 

\end{colorconjecture}
In pictures (see Fig. \ref{fig:flowers}), this means that the QR code output of any Whitehead double looks like a simple flower where the center constant is $6q$ and the coefficient of any petal is $-q$. If $q$ is positive, we see the red-centered flower; $q$ negative produces the blue-centered flower. 

Using information from KnotInfo \cite{KnotInfo} and The Knot Atlas \cite{TKA}, we verified the conjecture for all prime knots with fewer than 13 crossings (a total of 2977 knots). For each knot $K$, we computed the Whitehead double of $K$ using the algorithm described in Subsection \ref{codepart} and coded in \ref{appendix}. We then calculated the second derivative of the Alexander polynomial at $T=1$ and compared it to the result of $\theta (WD(+,K))$ divided by $\left( 6-T_1-T_2-T_1T_2-\frac{1}{T_1} - \frac{1}{T_2} -\frac{1}{T_1T_2}\right)$. Note that Whitehead doubles of $12$ crossing knots have at least $50$ crossings, further supporting that $\theta$ is genuinely computable. 

One important consequence shows that $\Theta$ may not be the best invariant for studying satellite knots. The conjecture implies that iterations of Whitehead doubles are zero, where other invariants can distinguish them.

\subsection{$t$-twisted Whitehead Double PD Algorithm}\label{codepart}
To efficiently compute $\theta$ for many examples of Whitehead doubles, it was helpful to write a code that took input as a knot and output the Whitehead double of the knot. To capture the crossings and orientation of the knot, we use planar diagram (PD) notation. The authors  acknowledge an independently and concurrently developed object-oriented code for use with SnapPy that generates Whitehead doubles written by Shana Li \cite{Li25}. In contrast, our code is purely functional and uses lists, making it easy to translate to other languages. Our intent is to share and explain the algorithm for future work on $t$-twisted Whitehead doubles. 

\begin{figure}
\centering
\begin{tikzpicture}
\begin{scope}[scale=0.8]
    \draw[line width=2pt, <-] (-1,1) -- (1,-1);
    \draw[line width=2pt] (-1,-1) -- (-0.1,-0.1);
    \draw[line width=2pt, <-] (1,1) -- (0.1,0.1);
    \node[below right] at (-0.75,-0.75) {$a$};
    \node[above right] at (0.75,-0.75) {$b$};
    \node[above left] at (0.75,0.75) {$c$};
    \node[below left] at (-0.75,0.75) {$d$};

    \draw[->] (2,0) -- (6,0);

    \begin{scope}[xshift=8cm,scale=2]
    \draw[line width=2pt,<-] (-1,1) -- (2,-2);
    \draw[line width=2pt,<-] (3,-1) -- (0,2);
    \draw[line width=2pt] (-0.1,-0.1) -- (-1,-1);
    \draw[line width=2pt] (0.9,0.9) -- (0.1,0.1);
    \draw[line width=2pt,<-] (2,2) -- (1.1,1.1);
    \draw[line width=2pt] (2.1,0.1) -- (3,1);
    \draw[line width=2pt] (1.9,-0.1) -- (1.1,-0.9);
    \draw[line width=2pt,->] (0.9,-1.1) -- (0,-2);
    \end{scope}

    \begin{scope}[xshift=8cm]
    \node[below right] at (-0.75-1,-0.75-1) {$w$};
    \node[below right] at (-0.75+1.8,-0.75+2) {$w+1$};
    \node[below right] at (-0.75+3.8,-0.75+4) {$w+2$};
    \node[above right] at (-0.75-1,0.75+1) {$x+2$};
    \node[above right] at (0.75,-0.9) {$x+1$};
    \node[above right] at (0.75+2.7,-0.75-2.7) {$x$};
    \node[above right] at (0.75,-0.75+4) {$z-2$};
    \node[above right] at (0.75+2,-0.75+2) {$z-1$};
    \node[above right] at (0.75+4.5,-0.75-0.5) {$z$};
    \node[below right] at (0.75,0.75-4) {$y$};
    \node[below right] at (0.75+2,0.75-2) {$y-1$};
    \node[below right] at (0.75+4,0.75) {$y-2$};
    \end{scope}

    \begin{scope}[xshift=4cm]
        \node at (0,1.5) {\scriptsize $w:=2a$};
        \node at (0,1.1) {\scriptsize $x:=2b$};
        \node at (0,0.7) {\scriptsize $y:=2a+(2n-a)\cdot4+2$};
        \node at (0,0.3) {\scriptsize $z:=2b+(2n-b)\cdot4+2$};
    \end{scope}
\end{scope}
\end{tikzpicture}
\caption{One crossing turns into four crossings}
\label{process-crossings}
\end{figure}
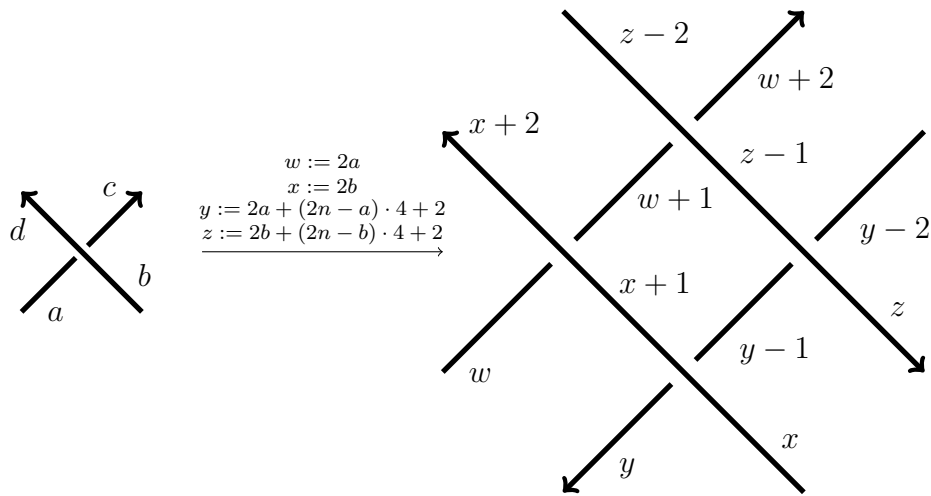

Planar diagram (PD) notation is a list of 4-tuples where each 4-tuple represents a crossing. First, edges are labeled $1$ to $2n$ following orientation along a knot diagram with $n$ crossings. The crossing 4-tuple is listed with its incoming lower edge first followed by edges that appear counterclockwise about the crossing \cite{KnotInfo}. For example, the crossing on the left of Figure \ref{process-crossings} is represented by $[a,b,c,d]$. Note that $b-d\equiv1\pmod{2n}$ for a positive crossing and $b-d\equiv-1\pmod{2n}$ for a negative crossing (as shown). 

Throughout the code, we use congruence modulo $2n$ where $0$ is redefined as $2n$ by defining safe\_mod in Python (\ref{python}) and SafeMod in Mathematica (\ref{mathematica}). We omit this for simplicity in this section. 

To create a diagram for the Whitehead double of $K$ starting from a diagram of $K$ with $n$ crossings, one can draw a parallel copy of $K$, add full twists to account for the writhe of $K$, and close with a positive or negative clasp. If the number of twists is greater than or less than the writhe by $t$ full twists, the resulting knot is called the $t$-twisted Whitehead double. \\\\
\noindent We now begin a walkthrough of the algorithms coded in \ref{appendix}.\\\\
\textbf{Step 1: The Set Up}\\
For input, we require a positive or negative sign for the clasp and PD notation for the knot $K$. For working in Mathematica, the KnotTheory` package \cite{TKA} allows one to simply provide the name (ex: Knot[3,1]) to call the PD notation for the knot. The number of twists is optional input with default set to zero twists to provide the untwisted Whitehead double. 

We capture the number of crossings in the diagram for $K$ using the length of the PD list. We calculate the writhe $Wr$ of $K$ and use it to calculate the number of edges in the Whitehead double 
\begin{equation}
    N = 8n + 4 + 4|t-Wr|
\end{equation}
where $8n$ edges come from each crossing becoming four new crossings, $+4$ edges make up the clasp, and $4|t-Wr|$ edges make up the full twists.

\begin{figure}[h]
    \centering
    \begin{tikzpicture}
       \node at (0,0) {\includegraphics[width=0.3\linewidth]{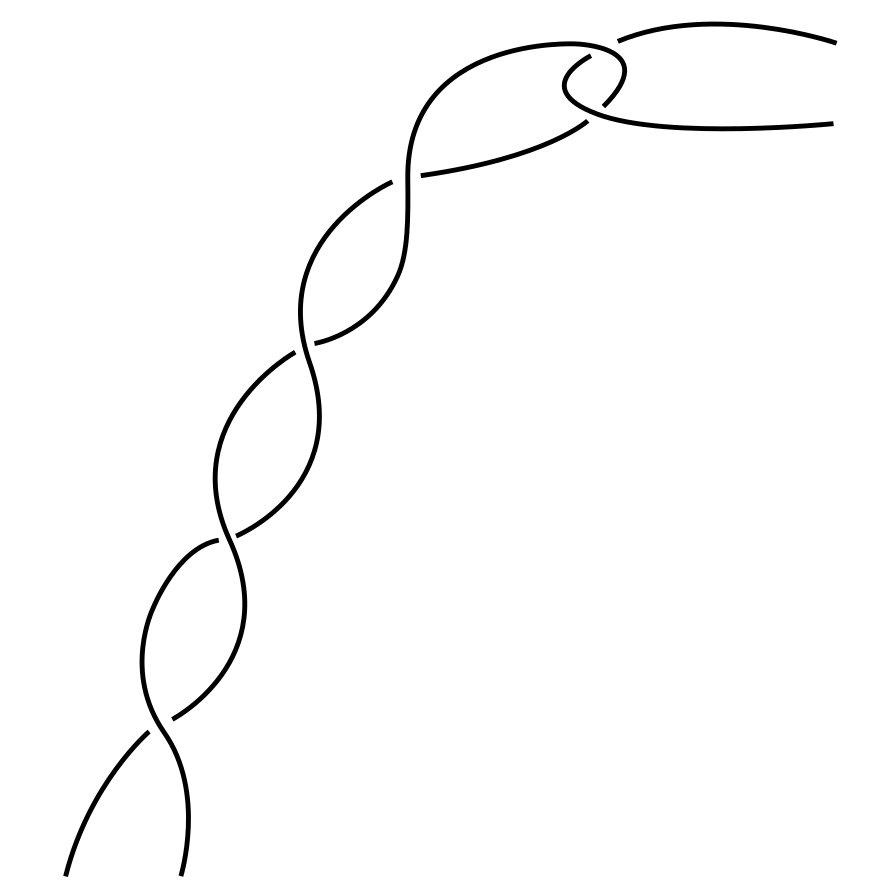}};
       \node at (0.5,-1) {two full twists};
       \draw[->] (-0.7,2) -- (-0.3,1.8);
       \node at (-1.1,2) {end};
       \draw[->] (1.3, 0.8) -- (1.3,1.1);
       \node at (1.3,0.6) {turnaround};
    \end{tikzpicture}
    \caption{The turnaround and end of the Whitehead double construction.} 
    \label{fig:turnaround}
\end{figure}

\textbf{Step 2: One Crossing to Four Crossings}\\
We label edges of four new crossings in terms of the lower edges of a crossing in $K$ as in Fig.~\ref{process-crossings}. We label our first edge by doubling $a$ ($w:=2a$). To get to the edge parallel to $w$, we must travel through the knot to the clasp ($(2n-a)\cdot 2$ edges), turn around ($+2$ edges), and travel back through the knot again ($(2n-a)\cdot 2$ edges). This parallel edge becomes $y:=2a+(2n-a)\cdot4+2$. We compute new edges dependent on $b$ similarly. We define special cases for when $a=2n$ or $b=2n$, which require traveling through the full twists to the end of the knot, instead of the turnaround. See Fig. \ref{fig:turnaround}.

\noindent \textbf{Step 3: Create Clasp}\\
We have two cases for positive or negative clasp. The clasp consists of the edges associated with the turnaround (edge $4n$) and the edges at the end of the twist section (edge $N - 1 - 2|t-Wr|$=: clasp\_center). The starting edge from Step 2 and the clasp divide the knot into full twists on one side and four-crossing sections on the other side. This placement simplifies the code considerably.\\\\
\textbf{Step 4: Create the Full Twists}\\
We have two cases for when $t-Wr$ is positive or negative. A for loop creates pairs of crossings (where each pair is one full twist) starting from the clasp\_center edge. The last twist aligns perfectly with the first edges $w$ and $y$ labeled in Step 2.

\appendix

\section{Associated Coded Functions}\label{appendix}

Both the Mathematica and Python code files are supplementary to this document. If not available as attachments, the essential functions appear in full below. However, the full files include examples and a graphical user interface in Python. 
\subsection{Python}\label{python}
\begin{verbatim}
def safe_mod(a,N):
    return ((a - 1) % N) + 1

def compute_whitehead_double(pd_code, clasp, twists):
    # takes inputs for any knot K, sign of clasp, 
    # and optional input integer t
    # to output the t-twisted *-clasped Whitehead double of K
    n = len(pd_code)
    t = twists
    WDPD = []

    # Find the writhe
    writhe = 0
    for crossing in pd_code:
        a,b,c,d = crossing
        if(b-d)%(2*n) == (-1)%(2*n):
            writhe -= 1
        else:
            writhe += 1

    # Process each crossing into four new crossings
    for crossing in pd_code:
        a,b,c,d = crossing
        # N is the number of edges in WD
        N = 8*n+4+4*(abs(writhe-t))
        # New edges w and y come from the original edge a
        # Similarly, x and z come from b
        w = 2*a
        x = 2*b
        y = 2*a + (2*n-a)*4+2
        z = 2*b + (2*n-b)*4+2
        if a == 2 * n:
            w = N
            y = N - 2 - 4 * abs(writhe - t)
        if b == 2 * n and (b - d) % (2 * n) == (-1) % (2 * n):
            x = N
            z = N - 2 - 4 * abs(writhe - t)

        if (b - d) % (2 * n) == (-1) % (2 * n):
            # Negative crossing
            WDPD.append([
                w,
                safe_mod(x + 1, N),
                safe_mod(w + 1, N),
                (x + 2) % N])
            WDPD.append([y - 1,safe_mod(x + 1, N), y, x])
            WDPD.append([y - 2, z - 1, y - 1, z])
            WDPD.append([
                safe_mod(w + 1, N),
                z - 1,
                safe_mod(w + 2, N),
                z - 2])
        else:
            # Positive crossing
            WDPD.append([
                w,
                z + 1,
                safe_mod(w + 1, N),
                safe_mod(z + 2, N)])
            WDPD.append([y - 1, safe_mod(z + 1, N), y, z])
            WDPD.append([y - 2, x - 1, y - 1, x])
            WDPD.append([
                safe_mod(w + 1, N),
                x - 1,
                safe_mod(w + 2, N),
                safe_mod(x - 2, N)])

    # Create clasp
    if clasp == '-':
        WDPD.append([
            4 * n,
            N - 1 - 2 * abs(writhe - t),
            4 * n + 1,
            N - 2 * abs(writhe - t)])
        WDPD.append([
            N - 2 - 2 * abs(writhe - t),
            4 * n + 1,
            N - 1 - 2 * abs(writhe - t),
            4 * n + 2])
    else: 
        WDPD.append([
            N - 1 - 2 * abs(writhe - t),
            4 * n + 1,
            N - 2 * abs(writhe - t),
            4 * n])
        WDPD.append([
            4 * n + 1,
            N - 1 - 2 * abs(writhe - t),
            4 * n + 2,
            N - 2 - 2 * abs(writhe - t)])
    clasp_center = N - 1 - 2 * abs(writhe - t)

    # First loop: writhe - t positive twists
    for i in range(1, writhe - t + 1):
        WDPD.append([
            clasp_center - 2*i - 1,
            clasp_center + 2*i + 1,
            clasp_center - 2*i,
            clasp_center + 2*i])
        WDPD.append([
            clasp_center + 2*i - 1,
            clasp_center - 2*i + 1,
            clasp_center + 2*i,
            clasp_center - 2*i])

    # Second loop: -writhe + t (negative) twists
    for i in range(1, -writhe + t + 1):
        WDPD.append([
            clasp_center - 2*i,
            clasp_center + 2*i - 1,
            clasp_center - 2*i + 1,
            clasp_center + 2*i])
        WDPD.append([
            clasp_center + 2*i,
            clasp_center - 2*i - 1,
            clasp_center + 2*i + 1,
            clasp_center - 2*i])

    return WDPD
\end{verbatim}

\subsection{Mathematica}\label{mathematica}
The following code appears in a supplementary file to this note or can be copied from below into Mathematica. The file contains examples using Bar-Natan's and van der Veen's $\Theta$, and the Mathematica notebook for computing $\Theta$ can be found at \cite{Ba25}. 
\begin{verbatim}
SafeMod[a_, n_] := Mod[a - 1, n] + 1

(*WD takes inputs sign s (1 or -1), any knot K, \
and optional input integer t to output the t-twisted, \
s-clasped Whitehead double of the knot K*)
WD[s_, K_, t_ : 0] := 
 Module[{WDPD = {}, writhe = 0, a, b, c, d, n, N, w, x, y, z},
  (*Number of crossings in K*)
  n = Length[PD[K]];
  
  (*Find the writhe*)
  Do[
   {a, b, c, d} = List @@ crossing;
   
   If[Mod[b - d, 2 n] == Mod[-1, 2 n], writhe -= 1, writhe += 1];
   ,
   {crossing, List @@ PD[K]}
   ];
  
  (*Process each crossing into four new crossings*)
  Do[
   {a, b, c, d} = List @@ crossing;
   
   (*Number of edges in WD*)
   N = 8 n + 4 + 4 (Abs[writhe - t]);
   (*New edges w and y come from the original edge a*)
   (*Similarly, x and z come from b, 
   this can be seen by overlaying the WD diagram on \
   top of the original diagram*)
   w = 2 a; x = 2 b; y = 2 a + (2 n - a) 4 + 2; 
   z = 2 b + (2 n - b) 4 + 2;
   If[a == 2 n, (w = N; y = N - 2 - 4 (Abs[writhe - t]))];
   If[(b == 2 n && Mod[b - d, 2 n] == Mod[-1, 2 n]),
    (x = N; z = N - 2 - 4 (Abs[writhe - t]))];
   
   (*Negative crossing*)
   If[Mod[b - d, 2 n] == Mod[-1, 2 n],
    AppendTo[WDPD, X[w, SafeMod[x + 1, N],
      SafeMod[w + 1, N], Mod[x + 2, N]]];
    AppendTo[WDPD, X[y - 1, SafeMod[x + 1, N], y, x]];
    AppendTo[WDPD, X[y - 2, z - 1, y - 1, z]];
    AppendTo[WDPD, 
     X[SafeMod[w + 1, N], z - 1, SafeMod[w + 2, N], z - 2]];
    ,
    (*Postive crossing*)
    AppendTo[WDPD, X[w, z + 1, SafeMod[w + 1, N], 
      SafeMod[z + 2, N]]];
    AppendTo[WDPD, X[y - 1, SafeMod[z + 1, N], y, z]];
    AppendTo[WDPD, X[y - 2, x - 1, y - 1, x]];
    AppendTo[WDPD, 
     X[SafeMod[w + 1, N], x - 1, SafeMod[w + 2, N],
     SafeMod[x - 2, N]]];
    ];
   ,
   {crossing, List @@ PD[K]}
   ];
  
  (*Create negative(s=-1) or positive(s=else) clasp*)
  If[s == -1,
   AppendTo[WDPD, 
    X[4 n, N - 1 - 2 (Abs[writhe - t]), 4 n + 1, 
     N - 2 (Abs[writhe - t])]];
   AppendTo[WDPD, 
    X[N - 2 - 2 (Abs[writhe - t]), 4 n + 1, 
     N - 1 - 2 (Abs[writhe - t]), 4 n + 2]];
   ,
   AppendTo[WDPD, 
    X[N - 1 - 2 (Abs[writhe - t]), 4 n + 1, 
     N - 2 (Abs[writhe - t]), 4 n]];
   AppendTo[WDPD, 
    X[4 n + 1, N - 1 - 2 (Abs[writhe - t]), 4 n + 2, 
     N - 2 - 2 (Abs[writhe - t])]];
   ];
  
  (*Create full-twists*)
  claspCenter = N - 1 - 2 (Abs[writhe - t]);
  For[i = 1, i <= writhe - t, i++,
   AppendTo[WDPD, 
    X[claspCenter - 2 i - 1, claspCenter + 2 i + 1, 
     claspCenter - 2 i, claspCenter + 2 i]];
   AppendTo[WDPD, 
    X[claspCenter + 2 i - 1, claspCenter - 2 i + 1, 
     claspCenter + 2 i, claspCenter - 2 i]]];
  For[i = 1, i <= -writhe + t, i++,
   AppendTo[WDPD, 
    X[claspCenter - 2 i, claspCenter + 2 i - 1, 
     claspCenter - 2 i + 1, claspCenter + 2 i]];
   AppendTo[WDPD, 
    X[claspCenter + 2 i, claspCenter - 2 i - 1, 
     claspCenter + 2 i + 1, claspCenter - 2 i]]];
  
  PD @@ WDPD
  ]
\end{verbatim}

\section*{Acknowledgments}
 The authors thank Dr. Matthew Hedden and Dr. Efstratia Kalfagianni for their advice and feedback during the project, Dr. Dror Bar-Natan and Dr. Roland van der Veen for their correspondence, and Dr. Christopher St. Clair for testing the Whitehead double Python GUI for user-friendliness. The second author thanks PhD advisors Dr. Matthew Hedden and Dr. Efstratia Kalfagianni. The second author is currently supported by funding from the NSF RTG Grant DMS-2135960.

\printbibliography

\end{document}